\documentstyle{amsart}

\numberwithin{equation}{section}

\newtheorem{theorem}{Theorem}[section]
\newtheorem{theorema}{Theorem}
\newtheorem{lemma}[theorem]{Lemma}
\newtheorem{corollary}[theorem]{Corollary}
\newtheorem{proposition}[theorem]{Proposition}
\newtheorem{corol}[theorema]{Corollary}

\theoremstyle{remark}
\newtheorem{remark}[theorem]{Remark}   
\newtheorem{remarks}{Remarks}

\newtheorem{definition}{Definition}
\newtheorem{ack}{Acknowledgment}

\begin{document}

\title[Complete hyperk\"ahler $4n$-manifolds]
{Complete hyperk\"ahler $4n$-manifolds with a local tri-Hamiltonian ${\Bbb R}^n$-action}
\author{Roger Bielawski}


\address{ Department of Mathematics and Statistics\\University of Edinburgh\\King's Buildings\\Edinburgh EH9 3JZ}

\email{rogerb@@maths.ed.ac.uk}

\begin{abstract} We classify those manifolds mentioned in the title which have finite topological type. Namely we show that any such connected $M^{4n}$ is isomorphic to a hyperk\"ahler quotient of a flat quaternionic vector space ${\Bbb H}^d$ by an abelian group.  We also show that a compact connected and simply connected $3$-Sasakian manifold of dimension $4n-1$ whose isometry group has rank $n+1$ is isometric to a $3$-Sasakian quotient of a sphere by a torus. As a corollary, a compact connected quaternion-K\"ahler $4n$-manifold with positive scalar curvature and  isometry group of rank $n+1$ is isometric to ${\Bbb H}P^{n}$ or $Gr_2({\Bbb C}^{n+2})$.
      
\end{abstract}

\maketitle

Hyperk\"ahler metrics in dimension $4n$ with $n$ commuting tri-holomorphic Killing vector fields form a large class of explicitly known Ricci-flat metrics.  The flat ${\Bbb H}^n$, the Taub-NUT, and the Eguchi-Hanson metric are of this form. So are the multi-Eguchi-Hanson metrics (a.k.a.  ALE-spaces of type $A_k$) and their Taub-NUT-like deformations, both due to Gibbons and Hawking. Higher dimensional examples include the Calabi metrics on $T^\ast {\Bbb C}P^n$ and the asymptotic metrics on moduli spaces of magnetic monopoles.
\par
A powerful method of constructing such metrics was given by Lindstr\"om and Ro\v{c}ek \cite{LR} (see also \cite{HKLR,PP}). It associates a hyperk\"ahler metric having a local tri-Hamiltonian (meaning Hamiltonian for all three K\"ahler forms) action of ${\Bbb R}^n$ to every real-valued function on an open subset of ${\Bbb R}^3\otimes {\Bbb R}^n$ which is {\em polyharmonic}, i.e. harmonic on any affine $3$-dimensional subspace of the form $a+{\Bbb R}^3\otimes {\Bbb R}v$, $a\in{\Bbb R}^3\otimes {\Bbb R}^n$, $v\in {\Bbb R}^n$. This construction gives locally {\em all} hyperk\"ahler metric with a local effective tri-Hamiltonian (hence isometric) action of ${\Bbb R}^n$. On the other hand, the global properties of the resulting metrics - in particular completeness - have not been investigated.
\par
There is another construction - that of hyperk\"ahler quotient - which, while in this case more restrictive, has the advantage that the global properties of the resulting manifold are accessible. In particular, a hyperk\"ahler quotient of a complete manifold by a compact Lie group is complete. A basic example is the hyperk\"ahler quotient of ${\Bbb H}\times (S^1\times {\Bbb R}^3)$ by the diagonal circle action which yields the Taub-NUT metric. More generally, we obtain complete hyperk\"ahler $4n$-manifolds with a tri-Hamiltonian (hence isometric) $T^n$-action by taking hyperk\"ahler quotients of quaternionic vector spaces by tori. This class of manifolds has been investigated in detail in \cite{BD}. They are closely related to projective toric varieties. In particular, they are diffeomorphic to a union of cotangent bundles of finitely many such varieties glued together.
\par
There is an interesting relation to symplectic geometry and Delzant's theorem \cite{De} which states that a compact symplectic $2n$-manifold with a Hamiltonian $T^n$-action is completely determined, as a Hamiltonian $T^n$-manifold, by the image of the moment map. Hyperk\"ahler $4n$-manifolds with a tri-Hamiltonian $T^n$-action can be thought of as the quaternionic analogue of that situation. Since now we cannot have compactness, we impose the next best thing: completeness of the metric.
\par
By the above-mentioned construction of Lindstr\"om and Ro\v{c}ek there are plenty of hyperk\"ahler manifolds with a local tri-Hamiltonian action of ${\Bbb R}^n$. We shall see that the situation is very different if we require the manifold to be complete (and Riemannian). On such a manifold any Killing vector field is complete and so an effective local tri-Hamiltonian action of ${\Bbb R}^n$ gives rise to an effective tri-Hamiltonian action of the connected abelian Lie group  $G={\Bbb R}^{p}\times T^{n-p}$.
We make the following definition.

\begin{definition} Let $M,M^\prime$ be two  hyperk\"ahler manifolds with tri-Hamiltonian actions of a Lie group $G$ and let $\mu,\mu^\prime$ be the chosen moment maps. We say that $M$ and $M^\prime$ are isomorphic as tri-Hamiltonian hyperk\"ahler $G$-manifolds, if there is a triholomorphic $G$-equivariant isometry $f:M\rightarrow M^\prime$ such that $\mu=\mu^\prime\circ f$.\end{definition}   

 The main result of this paper is:

\begin{theorema} Let $M^{4n}$ be a connected complete hyperk\"ahler manifold of finite topological type with an effective tri-Hamiltonian action of $G={\Bbb R}^{p}\times T^{n-p}$. Then 
\begin{itemize}
\item[(i)] If $M$ is simply connected and $p=0$, then $M$ is isomorphic, as a tri-Hamiltonian hyperk\"ahler $T^n$-manifold, to a hyperk\"ahler quotient of some flat ${\Bbb H}^{d}\times {\Bbb H}^m$, $m\leq n$, by $T^{d-n}\times {\Bbb R}^m$. 
\item[(i)] If $M$ is simply connected and $p>0$, then $M$ is isomorphic, as a tri-Hamiltonian hyperk\"ahler $G$-manifold, to the product of a flat ${\Bbb H}^p$ and a $4(n-p)$-dimensional manifold described in part (i). 
\item[(iii)] If $M$ is not simply connected, then $M$ is the product of a flat $(S^1\times {\Bbb R}^3\bigr)^l$, $1\leq l\leq n$, and a $4(n-l)$-dimensional manifold described in part (ii).\end{itemize} \label{1}\end{theorema}

\begin{remarks} 1. By a finite topological type we mean that the total Betti number is finite. There are examples \cite{AKN,Go} of complete $T^n$-invariant hyperk\"ahler manifolds not satisfying this condition.\\
2. In part (i), $T^{d-n}$ acts on ${\Bbb H}^d$ as a subtorus of the diagonal maximal torus $T^d$ in $Sp(d)$, while ${\Bbb R}^m$ acts on ${\Bbb H}^m\simeq \bigl({\Bbb R}^m\bigr)^4$ by translations on the first factor and on ${\Bbb H}^d$ via an injective linear map $\rho:{\Bbb R}^m\rightarrow \text{Lie}\,(T^d)={\Bbb R}^d$. In parts (ii) and (iii), the actions of ${\Bbb R}^p$ and $T^l$ are defined analogously, except that $\rho$ does not have to be injective.\\
3. Any {\em tri-symplectic} $T^n$-action on a simply-connected $M$ is tri-Hamiltonian.\\
4. In part (i), we can assume that the level set of the moment map used to obtain $M$ from ${\Bbb H}^{d+m}$ is smooth (this condition is always assumed in \cite{BD} and by no means follows from the smoothness of $M$ - see \cite{BGM2} for the $3$-Sasakian case).
 \end{remarks}

\begin{corol} Let $M$ be a simply connected $4$-dimensional complete hyperk\"ahler manifold  with a nontrivial tri-Hamiltonian vector field. If $b_2(M)=k>0$, then $M$ is isometric either to an ALE-space of type $A_k$ (i.e. a multi-Eguchi-Hanson space) or to its Taub-NUT-like deformation (i.e. to the hyperk\"ahler quotient  by ${\Bbb R}$ of the product of such a space with  ${\Bbb H}$). If $b_2(M)=0$, then $M$ is the flat ${\Bbb H}$. \label{2}\end{corol}
Here, if $b_2(M)>0$, we obtain an $S^1$-action, while if $b_2(M)=0$, the vector field gives rise either to the standard circle action or to translations in the first factor of ${\Bbb R}^4$.
\par

We emphasize that, for $n>1$, the $n$-dimensional  toral hyperk\"ahler quotients of  ${\Bbb H}^d$ are also quite well understood \cite{BD}. We know their integer cohomology, their homotopy type, generic complex structure and we have explicit formulas for the K\"ahler potential and the metric.

The previous corollary suggests the following definition.
\begin{definition} Let $M^{4n}$ be as in Theorem 1. A Taub-NUT deformation (of order $m$) of $M$ is the hyperk\"ahler quotient of $M\times {\Bbb H}^m$ by ${\Bbb R}^m$ where ${\Bbb R}^m$ acts on $M$ via an injective linear map $\rho:{\Bbb R}^m\rightarrow \text{Lie}\,(T^n)={\Bbb R}^n$. \end{definition}
Such a deformation $M^\prime$ is canonically $T^n$-equivariantly diffeomorphic to $M$ by a diffeomorphism $f$ which respects the hyperk\"ahler moment maps $\mu,\mu^\prime$, i.e. $\mu=\mu^\prime\circ f$.

We have the following hyperk\"ahler analogue of Delzant's theorem.
\begin{theorema} Complete connected hyperk\"ahler tri-Hamiltonian $T^n$-manifolds of finite topological type and dimension $4n$ are classified, up to Taub-NUT deformations, by arrangements of codimension $3$ affine subspaces in ${\Bbb R}^{3}\otimes {\Bbb R}^{n}$ of the form  
$$ H_k=\{(x^1,x^2,x^3)\in {\Bbb R}^{3}\otimes {\Bbb R}^{n};\, \langle x^i,u_k \rangle=\lambda_k^i,\enskip i=1,2,3\}$$
for some finite collection of vectors $u_k$ in ${\Bbb R}^n$ and scalars $\lambda_k^i$, $i=1,2,3$, such that, for any $p\in {\Bbb R}^{3}\otimes {\Bbb R}^{n}$, the set $\{u_k; p\in H_k\}$ is part of a ${\Bbb Z}$-basis of ${\Bbb Z}^n$.\label{3}\end{theorema}

This arrangement is the image of singular $T^n$-orbits by the hyperk\"ahler moment map. Taub-NUT deformations of different order are distinguished by the volume growth of the metric. In particular we have:
\begin{theorema}  Let $M$ be as in Theorem 1. If $M$ has Euclidean volume growth, then $M\simeq {\Bbb H}^p\times N$, where $N$ is isomorphic, as a tri-Hamiltonian hyperk\"ahler $T^{n-p}$-manifold, to a hyperk\"ahler quotient of a flat ${\Bbb H}^d$ by a $(d-n+p)$-dimensional subtorus $N$ of $Sp(d)$, i.e. to one of the toric hyperk\"ahler manifolds of \cite{BD}.\label{4}\end{theorema}
Remark 4 after Theorem 1 holds here as well.
\par
Sometimes one is interested in having an $SU(2)$  action which preserves the metric but not the complex structures. This is the situation for the Taub-NUT and for the metrics which are cones over $3$-Sasakian manifolds \cite{BGM,BD}. 

\begin{theorema} Let $M$ be as in Theorem 1 and suppose, in addition, that there is an effective isometric action of $SU(2)$ or $SO(3)$ rotating the complex structures and commuting with the action of $G$. Then $M$ is isometric to ${\Bbb H}^l\times (S^1\times {\Bbb R}^3)^{n-l}$, for some $l\leq n$, or its Taub-NUT deformation.\label{5}\end{theorema}

Our last two results deal with compact Einstein manifolds which admit some sort of quaternionic structure.
\par
We recall that a Riemannian manifold $(S,g)$ is $3$-Sasakian if the metric cone $({\Bbb R}_{>0}\times S, dr^2+r^2 g)$ is hyperk\"ahler. Such a manifold is Einstein of positive scalar curvature and examples of compact $3$-Sasakian $T^n$-manifolds were given and studied in \cite{BGM,BGMR,BD,BielBull}. These examples were obtained as $3$-Sasakian quotients (see \cite{BGM}) of spheres. We shall prove:

\begin{theorema} Let $N^{4n-1}$ be a connected simply connected compact $3$-Sasakian manifold admitting an effective $T^n$-action by $3$-Sasakian isometries (equivalently: the isometry group of $N$ has rank at least $n+1$). Then $N$ is isomorphic, as a $3$-Sasakian $T^n$-manifold, to a $3$-Sasakian quotient of a sphere by a torus.\label{6}\end{theorema}

In fact, the isomorphism can be required to preserve the $3$-Sasakian moment map for the $T^n$-action.
\par
A $3$-Sasakian manifold fibers over a quaternion-K\"ahler (i.e. having holonomy in $Sp(n-1)Sp(1)$) orbifold. Theorem \ref{6} in conjunction with Corollary 3 of \cite{BielBull} gives:

\begin{corol} A $4n$-dimensional compact connected quaternion-K\"ahler manifold with positive scalar curvature and isometry group of rank at least $n+1$ is isometric to ${\Bbb H}P^{n}$ or $Gr_2({\Bbb C}^{n+2})$ with their symmetric metrics.\label{7}\end{corol} \bigskip

The idea behind the proof of Theorem 1 is quite transparent. We shall first show (section \ref{red}) that it is sufficient to consider actions by tori. Then we observe that it is essentially sufficient to prove the result in dimension $4$, since, were there any new complete $T^n$-invariant hyperk\"ahler manifolds in higher dimensions, we would get new examples in dimension $4$ by taking hyperk\"ahler quotients by subtori. 
In dimension $4$, under the assumption of finiteness of the Betti numbers, the action of $S^1$ has finitely many fixed points $m_1,\dots, m_d$ and the hyperk\"ahler moment map induces a conformal immersion $\bar{\mu}:\bigl(M-\{m_1,\dots,m_d\}\bigr)/S^1\rightarrow {\Bbb R}^3$. Furthermore the local conformal factors are harmonic functions which implies that the metric on $\bigl(M-\{m_1,\dots,m_d\}\bigr)/S^1$ has nonnegative scalar curvature. We modify this metric so that it becomes a complete metric on $\bigl(M-\{m_1,\dots,m_d\}\bigr)/S^1$ with nonnegative scalar curvature. Then the results of Schoen and Yau \cite{SY} imply that $\bar{\mu}$ is injective and $\partial \mu(M)$ has Newtonian capacity zero (meaning that $\partial \mu(M)$ is removable for bounded harmonic functions). If $\mu(M)={\Bbb R}^3$, then we are done, since positive harmonic functions with finitely many isolated singularities at $\mu(m_1),\dots,\mu(m_d)$ are easily classified. To show that the boundary of $\mu(M)$ is indeed empty, we prove an estimate which implies that a metric, given by a harmonic conformal factor on a complement of a set $E$ of Newtonian capacity zero, is complete if and only if $E$ is empty.

\section{The Legendre transform and ${\Bbb R}^n$-invariant hyperk\"ahler metrics\label{LR}}

The fundamental idea of Lindstr\"om and Ro\v{c}ek \cite{LR} is that hyperk\"ahler manifolds of dimension $4n$ with a local effective tri-Hamiltonian action of ${\Bbb R}^n$ can be constructed from real-valued functions $F$ on an open subset $U$ of  ${\Bbb R}^3\otimes {\Bbb R}^n$ which are harmonic on $U\cap L$ for any $3$-dimensional affine subspace $L$ of the form ${\Bbb R}^3\otimes {\Bbb R}v$, $v\in {\Bbb R}^n$. If we identify ${\Bbb R}^3\otimes {\Bbb R}^n$ with ${\Bbb R}^n\times {\Bbb C}^{2n}$ with coordinates $x_i,z_i$, $i=1,\dots,n$, then Lindstr\"om and Ro\v{c}ek show that the Legendre transform with respect to the $x_i$ gives a K\"ahler potential $K$ of a hyperk\"ahler metric. The local holomorphic coordinates for this metric are given by the $z_i$ and some $u_i$ and we have
\begin{equation} K=F-2\sum_{i=1}^n(u_i+\bar{u}_i)x_i\label{K}\end{equation}
where the $x_i$ are determined by
\begin{equation} \frac{\partial F}{\partial x_i}=2(u_i+\bar{u}_i).\label{xi}\end{equation}
The vector fields $\partial/\partial y_i=\sqrt{-1}(\partial/\partial \bar{u}_i- \partial/\partial u_i)$ are triholomorphic isometries and the moment map for their action is simply $(x_i,\Re{z}_i,\Im{z}_i)$.
\par
An important observation of Pedersen and Poon \cite{PP} is that the metric has the form
\begin{equation} g=\sum_{i,j}\Bigl(\Phi_{ij}(dx_idx_j+dz_id\bar{z}_j)+ (\Phi^{-1})_{ij}(dy_i+A_i)(dy_j+A_j)\Bigr)\label{metric}\end{equation}
where $\Phi_{ij}=\frac{1}{4}F_{x_ix_j}$ and $A_j=\frac{\sqrt{-1}}{2}\sum_l ( F_{x_j\bar{z}_l}d\bar{z}_l - F_{x_jz_l}dz_l)$. In particular, the quotient metric on $U$ is given by
\begin{equation}\sum_{i,j}\Phi_{ij}(dx_idx_j+dz_id\bar{z}_j).\end{equation}
The $\Phi_{ij}$ are also polyharmonic functions and, in fact, the matrix $[\Phi_{ij}]$ determines (locally) the hyperk\"ahler and tri-Hamiltonian structure of the manifold. Indeed, if we have another function $F^\prime$ with $F^\prime_{x_ix_j}=4\Phi_{ij}$, then $G=F-F^\prime$ is of the form $\sum x_i g_i$, where the $g_i$ are polyharmonic functions of $z_j,\bar{z}_j$ only. Computing the difference of the two connection 1-forms $A$ and $A^\prime$, we see that it is given by 
$$\frac{\sqrt{-1}}{2}\sum_{i,j}\bigl(\frac{\partial g_i}{\partial\bar{z}_j}d\bar{z}_j- \frac{\partial g_i}{\partial z_j}dz_j\bigr).$$
Since the functions $g_i$ are polyharmonic, this form is closed and so it gives rise to a local equivariant isometry $\phi$. Moreover $\phi$ respects the K\"ahler potentials given by \eqref{K} and so the complex structure $I_1$. It also respects the coordinates $x_i,z_i$ (hence the moment map) and so, by the formula (2.8) in \cite{HKLR}, $\phi$ respects $I_2$ and $I_3$.\bigskip

The two basic examples of this construction are flat $S^1$-invariant metrics on $S^1\times {\Bbb R}^3$ and on ${\Bbb H}$. In the first case we have
\begin{equation} F(x,z,\bar{z})=2x^2- z\bar{z} \label{F1}\end{equation}
and consequently $\Phi\equiv 1$, while in the second case
\begin{equation} F(x,z,\bar{z})=x\ln (x+r)-r\label{F2}\end{equation}
where $r^2=x^2+z\bar{z}$. This time $\Phi=1/4r$. More general forms are given \cite{GM} and \cite{BD}. In the latter, the functions $F$ and the metrics for hyperk\"ahler quotients of flat vector spaces are computed. They are essentially obtained by taking linear combinations and compositions with linear maps of the solution \eqref{F2} (see also the proof of Theorem 1 in section \ref{proofs}). Our aim is to show that, in the case of a complete metric, the only other possibility is adding a linear combination of \eqref{F1}, which corresponds to a Taub-NUT deformation (see definition 2). 
\par
For a metric of the form \eqref{metric} taking hyperk\"ahler quotients by subtori is very simple. Indeed, the moment map equations are now linear (in $x_i,z_i$), and the hyperk\"ahler quotient corresponds to restricting the function $F$ to an appropriate affine subspace of ${\Bbb R}^3\otimes {\Bbb R}^n$. In fact, the requirement that $F$ be polyharmonic is a consequence of the fact that we must be able to take hyperk\"ahler quotients by any subtorus.
\par
An explanation of this construction in terms of twistors was given by Hitchin, Karlhede, Lindstr\"om and Ro\v{c}ek \cite{HKLR}. In particular, they have shown  that any hyperk\"ahler $4n$-manifold with a free tri-Hamiltonian  ${\Bbb R}^n$-action which extends to a ${\Bbb C}^n$-action  with respect to each complex structure and such that the moment map is surjective is given by the Legendre transform. In the next section we shall show that any hyperk\"ahler $4n$-manifold with a free local tri-Hamiltonian ${\Bbb R}^n$-action is {\em locally} given by the Legendre transform.

\section{Potentials for hyperk\"ahler metrics}

Let $X^{2n}$ be a K\"ahler manifold with the K\"ahler form $\psi$.  
It is well known that $\phi$ can be always written in a local holomorphic chart as
\begin{equation} \psi=-\frac{i}{2}\sum_{i,j} \frac{\partial^2K}{\partial z_i\partial \bar{z}_j}dz_i\wedge d\bar{z}_j\label{potential}\end{equation}
for a real-valued function $K$. 
\par
We now wish to show that a hyperk\"ahler structure is also locally determined by a single real-valued function $K$. Before proceeding we remark that the situation when there is a simultaneous K\"ahler potential for all three K\"ahler forms is very rigid \cite{Swann} - it is equivalent to the hyperk\"ahler manifold being locally a cone over a $3$-Sasakian manifold.
\par
Let $M^{4n}$ be an arbitrary hyperk\"ahler manifold and let $\omega_1, \omega_2,\omega_3$ be the three K\"ahler forms. Then $\omega=\omega_2+i\omega_3$ is holomorphic for the complex structure $I_1$. The Darboux theorem holds for complex symplectic forms and we can find a local $I_1$-holomorphic chart $u_i,z_i$, $i=1,\dots,n$ such that
\begin{equation} \omega=\sum_{i=1}^n du_i\wedge dz_i.\label{omega}\end{equation}
In this local chart $\omega_1$ can be written as in \eqref{potential}:
\begin{equation} \omega_1=-\frac{i}{2} \sum_{i,j}\left(K_{u_i\bar{u}_j}du_i\wedge d\bar{u}_j+ K_{u_i\bar{z}_j}du_i\wedge d\bar{z}_j+ K_{z_i\bar{u}_j}dz_i\wedge d\bar{u}_j +K_{z_i\bar{z}_j}dz_i\wedge d\bar{z}_j\right), \label{omega1}\end{equation}
for a real-valued function $K$. We see that the complex structure $I_2$ is given by:
\begin{align} I_2\left(\frac{\partial}{\partial u_i}\right) =\sum_{j=1}^n\left (K_{z_i\bar{u}_j}\frac{\partial}{\partial \bar{u}_j} + K_{z_i\bar{z}_j}\frac{\partial}{\partial \bar{z}_j}\right)\notag\\
I_2\left(\frac{\partial}{\partial z_i}\right) =\sum_{j=1}^n\left (-K_{u_i\bar{z}_j}\frac{\partial}{\partial \bar{z}_j} - K_{u_i\bar{u}_j}\frac{\partial}{\partial \bar{u}_j}\right).\label{J} \end{align}
Thus the condition $(I_2)^2=-1$ gives a system of nonlinear PDE's for $K$. This system is equivalent to the following condition:
\begin{equation}\begin{pmatrix} K_{u_i\bar{u}_j} & K_{u_i\bar{z}_j}\\ K_{z_i\bar{u}_j} & K_{z_i\bar{z}_j}\end{pmatrix} \in Sp(n,{\Bbb C}),\label{sp}\end{equation}
where the symplectic group is defined with respect to the form \eqref{omega}.
\par
 Conversely, suppose that in some local coordinate system $u_i,z_i$ we have a K\"ahler form $\omega_1$ given by a K\"ahler potential $K$ such that this system of PDE's is satisfied. Then, if we define $\omega_2+i\omega_3$ by the formula \eqref{omega}, we obtain a hyperhermitian structure. However $\omega_2$ and $\omega_3$ are closed, and so, by Lemma 4.1 in \cite{AH}, $I_2$ and $I_3=I_1I_2$ are integrable and we have locally a hyperk\"ahler structure. Therefore there is 1-1 correspondence between K\"ahler potentials satisfying the above system of PDE's and local hyperk\"ahler structures.
\par
For $n=1$ the condition \eqref{sp} reduces to the  Monge-Amp\`ere equation
\begin{equation} K_{u\bar{u}}K_{z\bar{z}}-K_{u\bar{z}}K_{z\bar{u}}=1.\label{MA} \end{equation}
 In higher dimensions, narrowing the situation from Ricci-flat K\"ahler to hyperk\"ahler is equivalent to replacing a single equation with high-order nonlinearities by a system of equations with quadratic nonlinearities. \medskip

We now go back to the main subject of this paper and we assume that there is a free local tri-Hamiltonian action of ${\Bbb R}^n$ on $M$. This action will extend to a local ${\Bbb C}^n$-action with respect to $I_1$ (or any other complex structure). Since there is a moment map for this ${\Bbb C}^n$-action, we can locally identify $M$ with a neighbourhood of $0$ in ${\Bbb C}^n\times {\Bbb C}^n$, where the first factor corresponds to the action and the second one is given by the moment map. Since we have now an equivariant holomorphic retraction onto ${\Bbb C}^n\times \{0\}$, the proof of the equivariant Darboux theorem, as given in \cite{GS}, goes through for $\omega=\omega_2+i\omega_3$. Therefore we identify locally $M$, as a complex ${\Bbb C}^n$-Hamiltonian manifold, with $T^\ast{\Bbb C}^n\simeq {\Bbb C}^n\times {\Bbb C}^n$ with coordinates $u_i,z_i$ in such a way that the $u_i$ correspond to the action of ${\Bbb C}^n$ and the $z_i$ give the moment map for this action. Furthermore $\omega$ has the form \eqref{omega}. The K\"ahler form $\omega_1$ is given by a K\"ahler potential $K$ which we can assume to be ${\Bbb R}^n$-invariant, i.e. independent of $u_i-\bar{u}_i$, $i=1,\dots,n$. Let $x_i$ denote the moment map for the action of $\sqrt{-1}(\bar{u}_i-u_i)$ with respect to $\omega_1$. From the formula \eqref{omega1} (using the fact that $K_{u_j}=K_{\bar{u}_j}$):
\begin{equation} x_i=-\frac{1}{2}K_{u_i}\label{xinew}.\end{equation}
We define a function $F$ by
\begin{equation} F=K+\sum_{i=1}^n 2x_i(u_i+\bar{u}_i).\label{Fnew}\end{equation}
In particular $F$ is independent of  $u_i-\bar{u}_i$, $i=1,\dots,n$.
Differentiating \eqref{Fnew} with respect to $u_i$, $i=1,\dots,n$,  yields the system of equations:
$$\sum_{j=1}^n\left(\frac{\partial F}{\partial x_j}-2u_j-2\bar{u}_j \right)\frac{\partial x_j}{\partial u_i}=0,\qquad i=1,\dots,n.$$
We claim that the matrix $[\partial x_j/\partial u_i]$ is everywhere nonsingular. Indeed suppose that at some point $(\underline{x},\underline{z})$ we have $\sum a_j\partial x_j/\partial u_i=0$, $i=1,\dots,n$, for some scalars $a_1,\dots,a_n$. The expression $\sum a_j x_j$ is the moment map for a  Hamiltonian vector field $Y_a$. If we restrict to the ${\Bbb C}^n$-orbit $O$ of  $(\underline{x},\underline{z})$  we conclude that the moment map for $Y_a$ has a critical point in the K\"ahler manifold $O$. This implies that $Y_a$ vanishes at  $(\underline{x},\underline{z})$  and so $a_1=\dots=a_n=0$. Therefore the above equations imply equations \eqref{xi}. Thus we have shown that the K\"ahler potential of a hyperk\"ahler metric with a free local tri-Hamiltonian ${\Bbb R}^n$-action can be written in the form \eqref{K} for a function $F$ satisfying \eqref{xi}. We wish to show that $F$ is polyharmonic. Let us take locally any hyperk\"ahler quotient of $M$ by an $(n-1)$-dimensional subgroup of ${\Bbb R}^n$. This corresponds to restricting $F$ to a $3$-dimensional affine subspace in ${\Bbb R}^3\otimes {\Bbb R}^n$ of the form $a+{\Bbb R}^3\otimes {\Bbb R}v$. Similarily, by identifying the hyperk\"ahler and complex-symplectic quotient, we restrict $K$ to an appropriate subspace. The restricted $K$ is still of the form \eqref{K} for the restricted $F$. Therefore to show that, in the general case, the function $F$ is polyharmonic it is enough to show that for $n=1$ such a function is harmonic. In this case, formula (3.143) of \cite{HKLR}, gives:
$$\bar{\partial} \partial K=\bar{\partial}(F_z)\wedge dz+2du\wedge \bar{\partial}(x).$$
This yields
$$ K_{z\bar{z}}=F_{z\bar{z}}+F_{zx}\frac{\partial x}{\partial \bar{z}}.$$
We also have, from \eqref{xinew},
$$ K_{u\bar{u}}=-2\frac{\partial{x}}{\partial \bar{u}},\enskip K_{u\bar{z}}=-2\frac{\partial{x}}{\partial \bar{z}}.$$
Now, from \eqref{xi},
$$ 0=\frac{\partial(u+\bar{u})}{\partial \bar{z}}= F_{xx}\frac{\partial x}{\partial \bar{z}}+F_{x\bar{z}}\enskip \text{
and}\enskip 
1=\frac{\partial(u+\bar{u})}{\partial \bar{u}}=  \frac{1}{2}F_{xx}\frac{\partial x}{\partial \bar{u}}.$$
Therefore
$$K_{u\bar{u}}=-4(F_{xx})^{-1},\enskip K_{u\bar{z}}=2(F_{xx})^{-1}F_{x\bar{z}}, \enskip  K_{z\bar{z}}=F_{z\bar{z}}-(F_{xx})^{-1}F_{x\bar{z}}F_{xz}$$
and we see that the equation \eqref{MA} is equivalent to $4F_{z\bar{z}}+F_{xx}=0$ which is the same as $\Delta F=0$. Thus we have proved
\begin{proposition} Let $M^{4n}$ be a hyperk\"ahler manifold with a free local tri-Hamiltonian action of ${\Bbb R}^n$. Then $M$ is locally given by the Legendre transform of Lindstr\"om and Ro\v{c}ek.\hfill$\Box$\label{Local} \end{proposition}  

\begin{remark} It follows from the above arguments that if $(M^{4n},I,\omega_1)$  is a K\"ahler manifold which also has a complex-symplectic form $\omega$ and a free local action of ${\Bbb R}^n$ which is Hamiltonian for both $\omega_1$ and $\omega$, then $M$ is locally given by the Legendre transform for some function $F$. This function is polyharmonic if and only if $M$ is hyperk\"ahler.\label{cs}\end{remark}

\section{Reduction to torus actions\label{red}}

In this section, we shall show that, in order to prove Theorem \ref{1}, it is enough to consider actions by tori. Suppose that we have a complete a hyperk\"ahler manifold $M^{4n}$ with an effective tri-Hamiltonian action of $G=T^{n-p}\times {\Bbb R}^{p}$, $p\leq n$. First we claim 
\begin{lemma} The group $G$ is closed in the full isometry group $I(M)$. \end{lemma}
\begin{pf} Consider the closure $\bar{G}$ of $G$ in $I(M)$. Let $\mu$ be the hyperk\"ahler moment map for $G$. Since $\mu$ is $G$-invariant, it is $\bar{G}$-invariant.  Let $m$ be a point at which $G$ acts freely and such that $s=\mu(m)$ is a regular value of $\mu$. Then $\mu^{-1}(s)$ is an $n$-dimensional submanifold of $M$ which is a union of $n$-dimensional orbits of $G$. Therefore  $O=Gm$ is closed in $\mu^{-1}(s)$ and so it is also an orbit of $\bar{G}$.  Thus $O=G$ and $O=\bar{G}/H$ where $H$ is closed in $\bar{G}$. Let $H_0$ be the identity component of $H$. Since $\bar{G}$ is abelian, there is a closed subgroup $F$ of $\bar{G}$ such that $\bar{G}=H_0\times F$, and so $O=F/E$ where $E=H\cap F$ is discrete. Since $O=G$ is a subgroup of $\bar{G}$, we have that $E=1$ and $G=F$ is closed in $\bar{G}$. \end{pf} 

\begin{corollary} The factor ${\Bbb R}^{p}$ in $G=T^{n-p}\times {\Bbb R}^{p}$ acts freely on $M$. \label{free}\end{corollary}
\begin{pf} This follows from the above lemma and the fact that isotropy groups of $I(M)$ are compact.\end{pf}

We shall now show
\begin{lemma} Any discrete subgroup of ${\Bbb R}^{p}$ acts properly discontinuously on $M$.\end{lemma}
\begin{pf} By Proposition 4.4 in \cite{KN}, it is enough to show that such a subgroup $L$ acts discontinuously, i.e. for any sequence $l_i$ of distinct elements of $L$ and any $x\in M$, $l_ix$ does not have a limit point in $M$.
\par
We shall prove several simple statements from which this claim will follow.
\newline
{\bf Claim 1).} If the action of $G$ is locally free at $x$, then $x$ does not belong to the closure of any $G$-orbit other than $Gx$. 
\par
 Indeed, if $x$ belongs to the closure of an orbit $Gy$, then $\mu(y)=\mu(x)$, but $Gx$ is an open neighbourhood of $x$ in $\mu^{-1}\bigl(\mu(x)\bigr)$, and so $Gy=Gx$.\newline
{\bf Claim 2).} If $l_ix$ converges to $y$ for some sequence $l_i$ of distinct elements of $L$, then $\dim Gy<n$.
\par
Indeed, the previous claim shows that $G$ cannot act locally freely at $y$.\newline
{\bf Claim 3).} If $l_ix\rightarrow y$, then $l_iy$ also has a limit point. 
\par
To show this we observe the following estimates, where $d$ denotes the distance on $M$  induced by the Riemannian metric:
\begin{multline*} d(l_iy,l_jy)\leq d(l_iy,l_il_ix)+d(l_il_ix,l_il_jx) + d(l_il_jx,l_jl_jx)+d(l_jl_jx,l_jy)\\
=d(y,l_ix)+2d(l_ix,l_jx)+d(l_jx,y)\leq 3d(y,l_ix)+3d(l_jx,y)\end{multline*}
where we have used the fact that $L$ is an abelian group of isometries. The claim follows now from completeness of the metric.\newline
{\bf Claim 4).} If $H$ is the stabilizer of a point $y\in M$, then the set $M^H$ of $H$-fixed points of $M$ is a closed hyperk\"ahler submanifold of $M$ of dimension $4n-4\dim H$ with a tri-Hamiltonian action of $T^{n-p}/H\times {\Bbb R}^p$.
\par
This is obvious, given that the tangent space to $M^H$ at $y$ is spanned by vectors $\{X_\rho,I_1X_\rho,I_2X_\rho,I_3X_\rho;\rho\in {\frak f}\}$ where $\text{Lie}(G)=\text{Lie}(H)\oplus{\frak f}$ and $X_\rho$ is the Killing vector field corresponding to $\rho$. We also use here Corollary \ref{free}.
\par  
 We can now show that $L$ acts discontinuously on $M$. Suppose, by contradiction, that there is a sequence $\{l_i\}$ of distinct elements of $L$ and a point $x\in M$ such that $l_ix$ converges to a point $y\in M$. Then, using claims 2), 3) and
4), we conclude that there is a complete $4s$-dimensional, $p\leq s< n$, hyperk\"ahler manifold $N$ with a {\sl free} tri-Hamiltonian action of $T^{s-p}\times {\Bbb R}^p$ and such that $l_iz$ has a limit point for some $z\in N$. This contradicts Claim 1).\end{pf} 

Therefore we can quotient $M$ by a maximal lattice $L$ in ${\Bbb R}^p$ and obtain a complete hyperk\"ahler $4n$-manifold with a tri-Hamiltonian action of $T^n$.

\section{Local structure of the moment map}

We shall prove some facts about Hamiltonian torus actions on hyperk\"ahler manifolds which we shall need in the sequel. We start with a simple
\begin{lemma} Let $M^{4n}$ be a hyperk\"ahler manifold equipped with an effective tri-Hamiltonian action of $T^m$. Then $m\leq n$.\label{triv}\end{lemma}
\begin{pf} At any point $p\in M$ let $V$ denote the subspace of the tangent space spanned by the vector fields generated by the action of $T^n$. It is simple to check, using the three moment maps, that the subspaces $V,IV,JV,KV$ are mutually orthogonal.\end{pf}

The following fact was proved in \cite{BD} for hyperk\"ahler manifolds which are toral quotients of flat quaternionic vector spaces (compare also \cite{De} for the symplectic case).
\begin{proposition} Let $M^{4n}$ be a hyperk\"ahler manifold equipped with an effective tri-Hamiltonian action of $T^n$ with a hyperk\"ahler moment map $\mu=(\mu_1,\mu_2,\mu_3):M\rightarrow {\Bbb R}^{3}\otimes {\Bbb R}^{n}$. Let $p\in M$ and let $V$ be a $T^n$-invariant neighbourhood of $p$, small enough that any singular orbit (i.e. one with nondiscrete stabilizer) intersecting $V$ contains $p$. Then:
\begin{itemize}
\item[(i)] The image $\mu(V^s)$ of the union $V^s$ of singular $T^n$-orbits  in $V$ is of the form $\mu(V)\cap \bigcup_{S} H_k$ where  
\begin{equation} H_k=\{(x^1,x^2,x^3)\in {\Bbb R}^{3}\otimes {\Bbb R}^{n};\, \langle x^i,u_k \rangle=\lambda_k^i,\enskip i=1,2,3\}\label{Hk}\end{equation}
for some countable collection $S$ of vectors $u_k$ in ${\Bbb R}^n$ and scalars $\lambda_k^i$, $i=1,2,3$;
\item[(ii)] the isotropy group $T_p$ of $p$ is the connected subtorus of $T^n$ whose Lie algebra is spanned by $\{u_k\in S; \mu(p)\in H_k\}$; 
\item[(iii)] After rescaling, the set $\{u_k\in S; \mu(p)\in H_k\}$ is part of a ${\Bbb Z}$-basis of ${\Bbb Z}^n$ and the isotropy representation of $T_p$ is determined by this set.\end{itemize}\label{moment}\end{proposition}
\begin{pf} Let $\rho\in {\Bbb R}^n$ and let $M_\rho$ be any component of the fixed point set of the 1-parameter subgroup of $T^n$ generated by $\rho$, i.e. the corresponding vector field $X_\rho$ vanishes on $M_\rho$. Then $M_\rho$ is a hyperk\"ahler submanifold of $M$. It follows that, for $p\in M_\rho$ and $v\in T_p M$, we have $\langle d\mu_i(v),\rho\rangle=0$, $i=1,2,3$. Therefore $\langle \mu_i,\rho\rangle$ is constant on $M_\rho$, $i=1,2,3$, which proves that $\mu(V^s)$ is as stated. This also proves that the Lie algebra of $T_p$ is spanned by $\{u_k\in S; \mu(p)\in H_k\}$. 
\par 
We shall now prove the statement (iii) and the connectedness of $T_p$. Let $\dim T_p=i$ and let $H_0$ be the identity component of $T_p$. There is an $(n-i)$-dimensional  subtorus $H_1$ of $T^n$ such that $T^n=H_0\times H_1$. Let $\Gamma=T_p/H_0=H_1\cap T_p$, so that the $T^n$-orbit of $p$ is isomorphic to $H_1/\Gamma$.  The fixed point set of $T_p$ is a totally geodesic hyperk\"ahler submanifold $X$ of $M$ and we can identify a $T^n$-invariant neighbourhood of $p$ with a $T^n$-invariant neighbourhood of $X$ in the normal bundle $N=T_XM/T_XX$. Since the torus $H_1/\Gamma$ acts freely on $X$, we can identify $X$ (via moment maps) with a neighbourhood of $H_1/\Gamma$ in  $H_1/\Gamma\times {\Bbb R}^{3n-3i}$. Since $M$ and $X$ both have a hypercomplex structure, so does the fiber of $N$, which we then identify with ${\Bbb H}^i$. The isotropy group $T_p$ must act on ${\Bbb H}^i$ via a homomorphism $\rho:T_p\rightarrow Sp(i)$. Since the action of $T^n$ is effective, $\rho$ must be injective and hence $\Gamma=1$. 
\par 
By the existence of normal forms for Hamiltonian actions (see section 41 in \cite{GS})  even the moment map $\mu_1$ can be identified with a symplectic moment map for the action of $T^n$ on $T^\ast{H_1}\times {\Bbb R}^{2i} \times {\Bbb C}^{2i}$ (the action is trivial on the middle factor). Let $\phi_1$ denote the moment map for the action of $T_p=H_0$ on ${\Bbb C}^{2i}$. If the weights of the representation are $\alpha_1,\dots,\alpha_i$, then we have (cf. \cite{BD}):
\begin{equation} \phi_1(z_1,\dots,z_i,w_1,\dots,w_i)=\frac{1}{2}\sum_{k=1}^i\bigl( |z_k|^2-|w_k|^2\bigr)\alpha_k\label{momentmap}.\end{equation}
Now it is easy to see that the vector $u_k$ is orthogonal to $\alpha_j$ for $i\neq k$. Since the representation is an isomorphism, the $\alpha_i$ form a ${\Bbb Z}$-basis of ${\Bbb Z}^i$, and therefore the $u_k$, after rescaling, also form a ${\Bbb Z}$-basis of their span. This finishes the proof.\end{pf}

Finally, we need: 
\begin{proposition} Let $M^{4n}$ be a hyperk\"ahler manifold equipped with an effective tri-Hamiltonian action of $T^n$. Then the hyperk\"ahler moment map $\mu=(\mu_1,\mu_2,\mu_3):M\rightarrow {\Bbb R}^{3n}$ induces a local homeomorphism from the orbit space $M/T^n$ into ${\Bbb R}^{3n}$. \label{open}\end{proposition}
\begin{pf}
  First we consider a $T^n$-invariant neighbourhood $U$ of a $T^n$-fixed point $m$,  which by the existence of normal forms we can identify
 with a neighbourhood of the origin in ${\Bbb C}^{2n}$ and the moment map $\mu_1$ with the standard moment map 
$$(z,w)=(z_1,\dots,z_n,w_1,\dots,w_n)\mapsto \frac{1}{2}\bigl(|z_1|^2-|w_1|^2,\dots, |z_n|^2-|w_n|^2\bigr).$$
 Furthermore, we can require that at the origin the complex structures $I_1,I_2,I_3$ coincide with the standard complex structures $I_1^0,I_2^0,I_3^0$ of ${\Bbb H}^n$. Let $\mu_1^0,\mu_2^0,\mu_3^0$ be the corresponding standard moment maps so that $\mu_1^0=\mu_1 ,(\mu_2^0+i\mu_3^0)(z,w)=(z_1w_1,\dots,z_nw_n)$. In a neighbourhood of $0$ we have an expansion $I_2=I_2^0+$({\sl terms of degree $\geq 1$ in $z_i,\bar{z_i},w_i,\bar{w_i}$}), and similarily for $I_3$. Now, for any vector $v$ tangent to $U$, $d\mu_2(v)=-d\mu_1(I_3v)$ and $d\mu_3(v)=d\mu_1(I_2v)$. Therefore 
\begin{equation}\mu_2+i\mu_3=\mu_2^0+i\mu_3^0+{\rm (}\text{\sl terms of degree $\geq 3$ in $z_i,\bar{z_i},w_i,\bar{w_i}$}{\rm )}.\label{expansion}\end{equation}
 The ring of (real) $T^n$-invariant polynomials on ${\Bbb C}^{2n}$ is generated by $y^i_1=|z_i|^2,y_2^i=|w_i|^2,y_3^i=\Re (z_iw_i), y_4^i=\Im(z_iw_i)$, $i=1,\dots, n$, with relations 
\begin{equation} y_1^i y_2^i=(y_3^i)^2+(y_4^i)^2.\label{relations}\end{equation}
In particular, as the invariant polynomials are homogeneous of degree $2$, the terms of degree $3$ in \eqref{expansion} vanish. We conclude that the moment map $\mu_2+i\mu_3$ can be written, on a sufficiently small neighbourhood of $0$, as 
$$(\mu_2+i\mu_3)(z,w)=f\bigl(|z_1|^2,|w_1|^2,\Re (z_1w_1), \Im(z_1w_1), \dots, |z_n|^2,|w_n|^2,\Re (z_nw_n), \Im(z_nw_n)\bigr),$$
where $f:{\Bbb R}^{4n}\rightarrow {\Bbb C}^n$ is a $C^\infty$ function of the form
$$ f(y_1^1,y_2^1,y_3^1,y_4^1,\dots,y_1^n,y_2^n,y_3^n,y_4^n)=\bigl(y_3^1+iy_4^1, \dots, y_3^n+iy_4^n\bigr) +O(r^2)$$
with $r$ denoting the distance from the origin in ${\Bbb R}^{4n}$.   The image $S$ of the neighbourhood of $m$ in the orbit space $M/T^n$ is identified with a neighbourhood of $0$ in the variety given by the equations \eqref{relations} with additional constraints $y_1^i\geq 0$, $y_2^i\geq 0$. The map $\bar{\mu}:M/T^n\rightarrow {\Bbb R}^{3n}$ is of the form $\bar{\mu}^0+O(r^2)$, with 
$$\bar{\mu}^0(y_1^1,y_2^1,y_3^1,y_4^1,\dots,y_1^n,y_2^n,y_3^n,y_4^n)=\bigl((y_1^1-y_2^1)/2,y_3^1,y_4^1,\dots, (y_1^n-y_2^n)/2,y_3^n,y_4^n\bigr).$$ One easily checks that the function $g:{\Bbb R}^{3n}\mapsto {\Bbb R}^{4n}$ defined by
\begin{multline} g(p_1^1,p_2^1,p_3^1,\dots, p_1^n,p_2^n,p_3^n)= \left(p_1^1, -p_1^1+\sqrt{(p_1^1)^2+(p_2^1)^2+(p_3^1)^2}, p_2^1,p_3^1,\dots,\right. \\ \left. p_1^n, -p_1^n+\sqrt{(p_1^n)^2+(p_2^n)^2+(p_3^n)^2}, p_2^n,p_3^n\right)\label{inverse}\end{multline}
is the the inverse of $\bar{\mu}^0$ restricted to $S$. Since $g$ is locally Lipschitz, it follows that $|\bar{\mu}(u)-\bar{\mu}(v)|\geq C|u-v|$ for $u,v\in S$, providing that we restrict ourselves to a sufficiently small neighbourhood of $0$ in ${\Bbb R}^{4n}$. Thus we have shown that $\bar{\mu}$ is open and locally $1-1$ in a neighbourhood of a fixed point. Consider now an arbitrary point $m\in M$ and let $l=\dim \text{Stab}_{T^n}(m)$. From Proposition \ref{moment} $\mu(m)$ lies on exactly $l$ flats $H_k$. The vectors $u_k$ defining these $H_k$
are a part of ${\Bbb Z}$-basis and we can complete it to a full ${\Bbb Z}$-basis. Let $A$ be the subtorus of $T^n$ generated by the $n-l$ vectors used to complete the basis and let $\mu_A:M\rightarrow {\Bbb R}^{3n-3l}$ be the moment map for $A$. According to Proposition \ref{moment}, $A$ acts freely in a neighbourhood of $m$ and so any hyperk\"ahler quotient of $M$ by $A$ with the level set of the moment map lying close enough to $\mu_A(m)$ is a manifold (in a neighbourhood of the point induced by $m$). For any such hyperk\"ahler quotient $\mu_A^{-1}(\lambda)/A$, the point induced by $m$ is a fixed point for $T^n/A$ and so, by the previous argument, the moment map for the action of $T^n/A$ induces a locally $1-1$ mapping from $\mu_A^{-1}(\lambda)/T^n$ to ${\Bbb R}^{3l}$. Now, $\mu_A^{-1}(\lambda)$ can be viewed as $\mu^{-1}(V_\lambda)$, where $V_\lambda$ is a $3l$-dimensional affine subspace of ${\Bbb R}^3\otimes {\Bbb R}^n$ of the form ${\Bbb R}^3\otimes W_\lambda$ with $W_\lambda$ generated by the vectors $u_k$ such that $m$ lies on $H_k$ and passing through the point $(0,\lambda)$ in the chosen ${\Bbb Z}$-basis of ${\Bbb Z}^n$. It follows that $\bar{\mu}$ is $1-1$ in a neighbourhood of the orbit of $m$. Since the spaces are locally compact and $\bar{\mu}$ is continuous, $\bar{\mu}$ is a homeomorphism of this neighbourhood onto its image. The orbit space is a topological manifold (this follows from \eqref{inverse}) and so $\bar{\mu}$ is a local homeomorphism into ${\Bbb R}^{3n}$.   
\end{pf}

\section{Metrics with harmonic local conformal factors}

We wish to consider the following situation which generalizes the $4$-dimensional case of Theorem 1. In this section $Y$ is a complete Riemannian $(n+3)$-dimensional manifold with an effective isometric action of $T^n$ which is not free only at a finite number of points $y_1,\dots,y_d$. Furthermore there is a $T^n$-invariant map $\mu:Y\rightarrow {\Bbb R}^3$ which induces a local homeomorphism $\bar{\mu}:Y/T^n\rightarrow {\Bbb R}^3$ whose restriction to $N=(Y-\{y_1,\dots,y_d\})/T^n$ (with the quotient metric) is a conformal immersion  with all local conformal factors harmonic. In other words $N$ can be covered by open sets $V_\alpha$, on which $\bar{\mu}$ is $1-1$ and the quotient metric $g$ on $V_\alpha$ satisfies 
$$(\bar{\mu}^{-1})^\ast g=\Phi_\alpha g_e$$
 where $g_e$ is the flat metric on ${\Bbb R}^3$ and $\Phi_\alpha$ is a positive harmonic function on an open subset $\bar{\mu}(V_\alpha)$ of ${\Bbb R}^3$. We shall also assume, to simplify the proof, that the $\Phi_\alpha$ do not extend to any $y_i$, $i=1,\dots,d$.  We aim to show
\begin{proposition} Under the above assumptions, $\bar{\mu}$ is a homeomorphism onto ${\Bbb R}^3$ and $g$ is globally conformally flat with the conformal factor $\Phi(x)=\sum_{i=1}^d\frac{a_i}{|x-x_i|}+b$, where $x_i=\mu(y_i)$ and $a_i>0,b\geq 0$, $i=1,\dots,d$, are constants. \label{n=4}\end{proposition}
We remark that a completely analogous result holds (and has a simpler proof) if we replace ${\Bbb R}^3$ with any ${\Bbb R}^p$, $p>3$. \medskip

We start the proof by observing that, since $\bar{\mu}$ is a homeomorphism in a neighbourhood of each $y_i$,  we have one (positive) harmonic function $\Phi$ in a neighbourhood of each point $y_i$. We have the following simple lemma about such functions.

\begin{lemma} Let $\Phi$ be a nonnegative harmonic function in $U-\{0\}$ where $U$ is a neighbourhood of $0$ in ${\Bbb R}^3$. Then $\Phi=c/r+\phi$ where $c$ is a nonnegative constant, $r(x)=|x|$ and $\phi$ is harmonic in $U$.\label{est1}\end{lemma}
\begin{pf} Choose $\epsilon>0$ so that $B(0,2\epsilon)\subset U$. Let $x$ be a point of the sphere $S(0,\epsilon)$. The Harnack inequality applied to the ball $B(x,\epsilon-t/2)$ implies that $\Phi(tx)\leq Ct^{-2}\Phi(x)$ for all $t>0$, where $C$ depends only on $\epsilon$. Thus $\Phi(x)=O(|x|^{-2})$ in $B(0,\epsilon)$. Then it follows that $\Phi=\Phi_1+c/r+\phi$ where $c,r,\phi$ are as in the statement and $\Phi_1$ is a linear combination (with constant coefficients) of the first order partial derivatives of the fundamental solution $1/r$. Since $\Phi$ is nonnegative, $\Phi_1\equiv 0$.\end{pf}

We now deform the metric on $N$.
\begin{proposition} There exists on $N$ a complete metric of nonnegative scalar curvature for which $\bar{\mu}$ is a conformal immersion.\label{scalar}\end{proposition}
\begin{pf} A metric of the form $\Phi g_e$ on an open subset of ${\Bbb R}^3$ has nonnegative scalar curvature if and only if $\Phi^{1/4}$ is superharmonic (see, for example, \cite{SY2}, Chapter V).  Since, all powers $\Phi^\alpha$, $0<\alpha<1$, of a harmonic (or superharmonic) function $\Phi$ are superharmonic, we know that the quotient metric $g$ on $N$ has a nonnegative scalar curvature. This metric is incomplete, but the incompleteness occurs only at the points $y_1,\dots y_d$ (as the metric on $M$ is complete). Hence, if we can modify the metric $g$ so that it is complete at these points and still has a nonnegative scalar curvature, we will be done. Therefore we work in a neighbourhood $U$ of the origin in ${\Bbb R}^3$ where, by Lemma \ref{est1}, the metric $g$ is given by $(c/r+\phi)g_e$. By the last assumption made about $Y$, $c>0$.  For any function $\psi$ and any $\alpha\in (0,1)$ we have
\begin{equation} \Delta (\psi^\alpha)<0 \Longleftrightarrow \psi \Delta \psi + (\alpha-1)|\nabla \psi|^2 <0.\label{alpha}\end{equation}
Assume now that $\psi$ is of the form $f(r)+\phi$, where $\phi$ is harmonic in $U$ and $f$ approaches $+\infty$ as $r\rightarrow 0$. Then we estimate the right-hand side of \eqref{alpha} as
$$\psi\Delta \psi + (\alpha-1)|\nabla \psi|^2\leq (f+\phi)\Delta f + (\alpha-1)\bigl(\nabla f\cdot\nabla(f+2\phi)\bigr).$$
For any $\epsilon >0$ and $r$ small enough, $f+\phi\leq (1+\epsilon) f$ and $2|\phi_i|\leq \epsilon |f_i|$, $i=1,2,3$, where the subscript $i$ denotes the $i$-th partial derivative. Therefore the last estimate can be replaced by
$$\psi\Delta \psi + (\alpha-1)|\nabla \psi|^2\leq (1+\epsilon)\bigl(f\Delta f +(\alpha-1) \frac{1-\epsilon}{1+\epsilon} |\nabla f|^2\bigr).$$
In our case, $\alpha=1/4$, and if we choose $\epsilon$ so that $(\alpha-1) \frac{1-\epsilon}{1+\epsilon}=-1/2$, we see from \eqref{alpha} that if $f^{1/2}$ is superharmonic, then $\psi^{1/4}$ is superharmonic (in a small neighbourhood of the origin). 
\par
Our original conformal factor was of the form $c/r +\phi$. In order to make the metric complete we wish the growth to be $1/r^2$. Therefore we have to find a positive function $f(r)$ such that $f(r)=c/r$ for $r\geq \delta$, $f(r)\sim 1/r^2$ near $0$,  and $\sqrt{f}$ is superharmonic. Let $u=\sqrt{f}$. The condition $\Delta u\leq 0$ is equivalent, for a function $u=u(r)$, to 
$$\ddot{u}+\frac{2}{r}\dot{u}\leq 0,$$
and so, for $z=\frac{d}{dr}\ln(-{du}/{dr})$, it is equivalent to $z\leq 2/r$. Let $z(r)$ be any smooth (on $(0,\infty)$) decreasing function satisfying this inequality,  equal to $3/2r$ for $r\geq \delta$ and equal to $2/r$ for $r\leq \delta/2$. We obtain a function $u$ by choosing the two free constants (i.e. $u(\delta)$ and $\dot{u}(\delta)$) so that $u(r)$ is tangent to $\sqrt{c/r}$ at $r=\delta$. We define a new conformal factor by $\tilde{\Phi}=u^2+\phi$. It coincides with $\Phi=c/r +\phi$ for $r\geq \delta$ and is positive as $u(r)\geq \sqrt{c/r}$. Furthermore $\tilde{\Phi}^{1/4}$ is superharmonic and, therefore, $\tilde{\Phi}g_e$ has nonnegative scalar curvature. Finally, this metric is complete at $0$ as $\tilde{\Phi}$ has $(1/r^2)$-growth.\end{pf}

We now appeal to results of Schoen and Yau \cite{SY} (Propositions 4.2, 4.3, 4.4, and 4.4'; see also Theorem VI.3.5 in \cite{SY2}) which give us:
\begin{corollary} The map $\bar{\mu}:N\rightarrow {\Bbb R}^3$ is injective and the boundary of $\bar{\mu}(N)$ has Newtonian capacity zero.\hfill $\Box$\label{SYcor}\end{corollary}

We recall that a $G_{\delta}$ subset $E$ of ${\Bbb R}^n$ has Newtonian capacity zero (or is polar) precisely when the removable singularity theorem holds for $E$. Another equivalent condition is that there is a tempered positive measure $\mu_E$ in ${\Bbb R}^n$ such that the convolution of the Green function $G(x,y)=|x-y|^{2-n}$ with $\mu$ is infinite exactly on the set $E$. In particular the Hausdorff dimension of $E$ is at most $n-2$.\medskip 

Since the map $\bar{\mu}:Y/T^n \rightarrow {\Bbb R}^3$ is a local homeomorphism in a neighbourhood of each $y_i$ and since a set $E\cup \{x_1,\dots,x_d\}$, $x_i\not\in \bar{E}$, is polar if and only if $E$ is polar, the conclusion of the above corollary holds also for $Y/T^n$, instead of $N$. Therefore we can assume that $Y/T^n$ is an open subset $U$ of ${\Bbb R}^3$, and that the quotient metric on $U-\{x_1,\dots,x_d\}$, where $x_i=\mu(y_i)$, is of the form 
$$g=\left(\sum_{i=1}^d\frac{a_i}{|x-x_i|}+\phi\right)g_e$$
for positive constants $a_i$ and a function $\phi$, harmonic in $U$. 

We now wish to show that $\bar{\mu}:Y/T^n\rightarrow {\Bbb R}^3$ is onto. This will imply that $\phi$ is constant and will end the proof of Proposition \ref{n=4}. \par
We have:
\begin{lemma} Suppose that the metric on $Y$ is complete. Then $\phi$ is bounded below and for any $C^1$-curve $\gamma:[0,1]\rightarrow \bar{U}$, $\gamma\bigl([0,1)\bigr)\subset U$, $\gamma(1)\in \partial U$, we have
$$\int_0^1\bigl(\phi(\gamma(t))+C\bigr)^{1/2}|\gamma^\prime(t)|dt=+\infty$$
for any $C\geq 0$ such that $\phi+C\geq 0$.\label{vorth}\end{lemma}

In other words, for any $C\geq 0$ such that $\phi+C> 0$, $(\phi+C)g_e$ is a complete Riemannian metric on $U$.

\begin{pf} The function $\phi$ is bounded in a compact neighborhood $K_i$ of each of the $x_i$ and outside of $\bigcup \text{Int}\,K_i$ the function $f=\Phi-\phi$ is also bounded. Since $\phi+f>0$, $\phi$ must be bounded below.
\par
As in Proposition \ref{scalar}, the only incompleteness of the quotient metric $g$ occurs at the points $x_i$ and so a curve $\gamma$ with the assumed properties must have infinite length in $g$. We can assume that $\gamma$ avoids the points $x_i$. Let $\Psi$ be the restriction of $\sum a_k/|x-x_k|$ to $U$. This is finite on $\gamma$. For any $C\geq 0$ such that $\phi+C\geq 0$ we have
\begin{multline*} +\infty=\int_0^1\bigl((\Psi+\phi)(\gamma(t))\bigr)^{1/2}|\gamma^\prime(t)|dt\\
\leq \int_0^1\bigl(\Psi(\gamma(t)\bigr)^{1/2}|\gamma^\prime(t)|dt+ \int_0^1\bigl(\phi(\gamma(t))+C\bigr)^{1/2}|\gamma^\prime(t)|dt + \int_0^1 C^{1/2}|\gamma^\prime(t)|dt.\end{multline*}
The first and the last term are finite, and so the second one is infinite.\end{pf}

 Proposition \ref{n=4} will follow from the following fact (compare with Lemma \ref{est1}):
\begin{proposition} Let $E$ be a closed subset of Newtonian capacity zero in ${\Bbb R}^3$ and let $\Psi$ be a nonnegative harmonic function on ${\Bbb R}^3-E$. There is a constant $K\geq 0$ such that for all $x\in {\Bbb R}^3$,
\begin{equation}\Psi(x)\leq K+\frac{K}{\text{dist}_{S^3}(x,E)},\label{Psi-final}\end{equation}
where the distance is measured in the standard metric on $S^3$.\label{Phi-final}\end{proposition}
\begin{pf} To prove the estimate \eqref{Psi-final} for closed sets in ${\Bbb R}^3$ in the distance measured in $S^3$ is the same as to prove it for compact subsets of ${\Bbb R}^3$ in the Euclidean distance. Thus we can assume that $E$ is compact. Let now $\mu_E$ be the positive measure in ${\Bbb R}^3$ whose convolution with $|x-y|^{-1}$ gives the nonnegative superharmonic function  $\psi_E$ on ${\Bbb R}^3$ with $\psi_E^{-1}(+\infty)=E$. Then the function
$$f(x)=\begin{cases}\Psi(x)+ \psi_E(x) & \text{if $x\not\in E$}\\
+\infty & \text{if $x \in E$}\end{cases}$$
is nonnegative and superharmonic on ${\Bbb R}^3$. It follows (Corollary 2 of Theorem 1.23 in \cite{La}) that $f=c+G\ast \mu^\prime$, where $c$ is a nonnegative constant, $G(x,y)=|x-y|^{-1}$ and $\mu^\prime$ is a nonnegative measure in ${\Bbb R}^3$. Let $\nu=\mu^\prime -\mu_E$. Then $G\ast \nu$ coincides with $\Psi-c$ on ${\Bbb R}^3-E$. Since $\Psi$ solves the Laplace equation on ${\Bbb R}^3-E$ and $G\ast \nu$ is the solution of the Poisson equation $\Delta(G\ast \nu)=-\nu$ in the sense of distributions, we conclude that $\text{supp}\, \nu\subset E$. Let us write $\nu=\nu^+-\nu^-$ where $\nu^+,\nu^-$ are the positive and negative variations of $\nu$. We have
$$\Psi(x)-c=\int_E \frac{d\nu(y)}{|x-y|}\leq \int_E \frac{d\nu^+(y)}{|x-y|}+ \int_E \frac{d\nu^-(y)}{|x-y|}\leq \frac{\nu^+(E)+\nu^-(E)}{\text{dist}(x,E)}.$$
This concludes the proof.\end{pf}

We can now finish the proof of Proposition \ref{n=4}.\newline
{\em Proof of Proposition \ref{n=4}}. Suppose that $U=\mu(Y)\neq {\Bbb R}^{3n}$. We consider the metric $(\phi+C)g_e$ on $U$ for some large $C$. This metric is complete by Lemma \ref{vorth}.  Let $z\in {\Bbb R}^3$ be such that $\text{dist}(z,\partial U)=1$ and let $y\in \partial U$ be a point where this distance is achieved. Then, by Corollary \ref{SYcor}, the estimate \eqref{Psi-final} holds for all points $x$ of the segment $yz$. Therefore the length of this segment is finite in the metric $(\phi+C) g_e$ and we obtain a contradiction. Hence $U={\Bbb R}^3$ and $\phi$ must be a constant.\hfill $\Box$       \medskip

We finish the section by observing that, if we write  $(x,z)$ for the coordinates of ${\Bbb R}^3$ as in section \ref{LR}, then the function $\Phi$ of Proposition \ref{n=4} can be written as $\Phi=F_{xx}$ where (up to a term linear in $x$) $F$ is given by
\begin{equation} F(x,z,\bar{z})=\sum_{k=1}^da_k\bigl(s_k\ln(s_k+r_k)-r_k\bigr) +b\bigl(2x^2-z\bar{z}\bigr)\label{F3}\end{equation}
where $s_k=x-x_k$ and $r_k$ is the distance between $(x,z,\bar{z})$ and $\mu(y_k)=(x_k,z_k,\bar{z}_k)$.

\section{Proofs\label{proofs}}

We shall now prove all results stated in the introduction.

{\em Proof of Theorem 1.} Let $M$ satisfy the assumptions of Theorem 1 with $p=0$ (we showed in section \ref{red} that we can assume this). Let $H_k$, $k=1,\dots,d$, be the codimension $3$ affine subspaces in ${\Bbb R}^3\otimes {\Bbb R}^n$ given by Proposition \ref{moment}. Their number is finite because $M$ has finite topological type (this follows from the localization theorem \cite{Br}, if we notice that were there infinitely many $H_k$, we could find a circle in $T^n$ with infinitely many fixed points). Consider a small open subset $U$ of ${\Bbb R}^3\otimes {\Bbb R}^n$ not intersecting any $H_k$. In particular the $T^n$-action is free on $\mu^{-1}(U)$ and, if $U$ is small enough, the hyperk\"ahler moment map is a diffeomorphism between $\mu^{-1}(U)/T^n$ and $U$. Therefore, as in section \ref{LR}, there is a function $F=F(x_i,z_i,\bar{z}_i)$ defined on $U$ which completely determines a neighbourhood $V$ of a section of $\mu^{-1}(U)\rightarrow U$ as a hyperk\"ahler tri-Hamiltonian $T^n$-manifold. This function is unique up to irrelevant terms linear in the $x_i$.
\par
For any rational $3$-dimensional affine subspace $L=a+{\Bbb R}^3\otimes {\Bbb R}v$ of ${\Bbb R}^3\otimes {\Bbb R}^n$ such that $L\cap U\neq\emptyset$ we can consider the $(n-1)$-dimensional subtorus $N$ of $T^n$ whose Lie algebra is generated by vectors orthogonal to $v$. The level set of the moment map $\mu_N$ for $N$ is chosen to be $L$. According to Proposition \ref{moment}, as long as $L$ does not meet any point at which two of the $H_k$ intersect, $N$ will act locally freely on $\mu_N^{-1}(L)$ and so $Y=\mu_N^{-1}(L)$ is a complete Riemannian manifold which we claim satisfies all assumptions of the previous section. The map $\mu:Y\rightarrow L\simeq {\Bbb R}^3$ is simply the restriction of $\mu:M\rightarrow {\Bbb R}^3\otimes {\Bbb R}^n$. By Proposition \ref{open}, $\bar{\mu}$ is a local homeomorphism. The points $y_i$ map, under $\mu$, to the intersections of $L$ with the $H_k$. 
Therefore $Y=\mu_N^{-1}(L)$ satisfies all assumptions of the previous section and, by Proposition \ref{n=4}, we know, up to irrelevant linear factors, the restriction of $F$ to $L\cap U$ . It is given by  the formula \eqref{F3} where the points $(x^k,z^k)$ are intersection points of $L$ with the $H_k$.  Thus we know the restriction of $F$ to a generic rational $3$-dimensional subspace of $U$. This obviously determines $F$. We claim that on $U$, $F$ is given by the formula
\begin{equation}  F(x_i,z_i,\bar{z}_i)=\sum_{k=1}^d a_k \bigl(s_k\ln(s_k+r_k)-r_k\bigr) +\sum_{i,j}b_{ij}\bigl(2x_ix_j-z_i\bar{z}_j\bigr)\label{F4}\end{equation}
where $a_k$ are positive constants, $[b_{ij}]$ is a positive-definite constant matrix and $s_k$ and $r_k$ are defined by (cf. \eqref{Hk})
$$ s_k(x,z)=\langle x,u_k\rangle -\lambda_k^1,\enskip v_k(x,z)=\langle z,u_k\rangle -\lambda_k^2-\sqrt{-1}\lambda_k^2,\enskip r_k^2=s_k^2+v_k\bar{v}_k. $$
In other words $r_k$ is the distance of a point from $H_k$, and $|s_k|$ is the distance between the projections of a point and of $H_k$ onto the first factor in ${\Bbb R}^n\times {\Bbb R}^n\times {\Bbb R}^n\simeq {\Bbb R}^3\otimes {\Bbb R}^n$. Indeed, the second sum in \eqref{F4} is the most general form of a function which restricted to any (or generic rational) $3$-dimensional subspace $L$ gives us the second term of \eqref{F3}. For the first sum, notice that \eqref{F3} implies that $F$ must be of the form
$$\sum_{k=1}^p a_k \bigl(s_k\ln(s_k+r_k)-r_k\bigr) +\sum_{i,j}b_{ij}\bigl(2x_ix_j-z_i\bar{z}_j\bigr)$$
where $s_k$ are linear in the $x_i$ and $r_k^2$ are quadratic in the $x_i,\Re z_i, \Im z_i$. If, however, any of these was not of the form stated above, then for some $3$-dimensional subspace $L$ we would have obtained singular points of $Y_L$ different from those corresponding to the intersections of $L$ with the $H_k$. Therefore $F$ is of the form \eqref{F4}.

The first sum in \eqref{F4} describes the function $F$ of a hyperk\"ahler quotient $M^\prime$ of a flat ${\Bbb H}^d$ by a torus determined by the $H_k$. This follows from \cite{BD}, given that the $H_k$ satisfy the ${\Bbb Z}$-basis condition of Proposition \ref{moment}. The metric on ${\Bbb H}^d$ is not the standard one but each factor is rescaled by $a_k$. 
\par
Suppose that $[b_{ij}]$ has rank $m$. Since it is symmetric, we can write $[b_{ij}]=AA^T$ for an $n\times m$  matrix $A$ of rank $m$. This $A$ defines an embedding ${\Bbb R}^m\rightarrow {\Bbb R}^n$. The function \eqref{F4} is the the function $F$ of the Taub-NUT modification (see Definition 2) of order $m$ of $M^\prime$. Let us denote $M^\prime$ with this modified metric by $M^{\prime\prime}$. We conclude that $M$ and $M^{\prime\prime}$ are locally (on connected subsets mapping onto $U$) isomorphic as tri-Hamiltonian hyperk\"ahler $T^n$-manifolds. Let $\phi$ be this local isomorphism. In particular $\phi$ is an isometry and, as $M$ and $M^{\prime\prime}$ are Einstein, hence real-analytic, and complete, $\phi$ extends to a (unique) isometry between the universal covers
$\tilde{M}$ and  $\tilde{M}^{\prime\prime}$ (\cite{KN}, Corollary VI.6.4). If there are $n$ flats $H_k$ intersecting in a point, then, by \cite{BD}, $M^{\prime\prime}$ is simply connected and so we obtain an isometry between $M$ and $M^{\prime\prime}$. This isometry is locally a tri-Hamiltonian isomorphism and so it is such an isomorphism globally (as everything is real-analytic). This proves case (i) of Theorem 1.
\par
If there are at most $n-l$, $l>0$, flats $H_k$ intersecting at any given point, then $M^{\prime\prime}$ is isomorphic to a product of a simply connected $X$ and $(S^1\times{\Bbb R}^3)^l$. Hence the universal covers are isometric to $X\times ({\Bbb R}^4)^l$. Once more the isometry between $\tilde{M}$ and $X\times ({\Bbb R}^4)^l$ respects the tri-Hamiltonian structure for the action of $T^{n-l}\times {\Bbb R}^l$ (${\Bbb R}^l$ acts as in Remark 2 after Theorem \ref{1}). It follows that $\pi_1(M)$ is a subgroup of $T^{n-l}\times {\Bbb R}^l$. Since $T^{n-l}$ has a fixed point in $X$, any element of $\pi_1(M)$ has a nontrivial component in ${\Bbb R}^l$ and, hence, the projection onto the second factor of $T^{n-l}\times {\Bbb R}^l$ gives us an embedding  of $\pi_1(M)$ into ${\Bbb R}^l$. Since we do have a $T^n$-action on $M$, $\pi_1(M)$ must be a full lattice of ${\Bbb R}^l$ which acts on $X$ via a rational homomorphism $\lambda:{\Bbb R}^l\rightarrow {\Bbb R}^{n-l}\simeq \text{Lie}\,(T^{n-l})$. We can choose a basis of ${\Bbb R}^l$ so that $\pi_1(M)$ is identified with ${\Bbb Z}^l$ and the $i$-th generator of  ${\Bbb Z}^l$ acts on $X$ as $t_i\in T^{n-l}$. Then the map $X\times {\Bbb R}^{4l}$ to itself given by $(x,y)\mapsto (t_1^{-1}x,y)$ induces an equivariant isomorphism between $\bigl(X\times {\Bbb R}^{4l}\bigr)/{\Bbb Z}$ and $X\times (S^1\times{\Bbb R}^3) \times  {\Bbb R}^{4(l-1)}$. We can continue with successive generators of ${\Bbb Z}^l$ and conclude that $M$ is isomorphic to some $X\times (S^1\times {\Bbb R}^3)^l$. This proves case (ii) of Theorem 1.\hfill $\Box$

{\em Proof of Corollary \ref{2}.} This follows from the classification of $4$-dimensional quotients of vector spaces by tori, see for example \cite{BD}. \hfill $\Box$

{\em Proof of Theorem \ref{3}.} This follows from Theorem 1, Proposition \ref{moment} and the results of \cite{BD}. \hfill $\Box$

{\em Proof of Theorem \ref{4}.} We have to compute the volume growth of a hyperk\"ahler quotient of a vector space by a torus or its Taub-NUT deformation. The metric on a hyperk\"ahler quotient of a vector space by a torus is asymptotic to a cone metric on a $3$-Sasakian space and therefore it has Euclidean volume growth. Now suppose that $M$ is a Taub-NUT deformation of such a quotient of order $m$. If we diagonalize the matrix $[b_{ij}]$, then $\Phi_{ij}=F_{x_ix_j}$ is of the form $\lambda_{ij}+\sum_{k=1}^d a_{ijk}/r_k$ with $\lambda_{ij}=0$ if $i\neq j$ or $i>m$. The volume growth is then comparable to the volume growth of $(S^1\times {\Bbb R}^3)^m\times (\text{\sl metric cone of dimension $4n-4m$})$, which is $4n-m$. \hfill $\Box$     

{\em Proof of Theorem \ref{5}.} Such an action of $SU(2)$ or $SO(3)$ induces the standard $SO(3)$-action on the first factor of ${\Bbb R}^3\otimes {\Bbb R}^n$. The affine subspaces $H_k$ must be preserved by this action, and so they all pass through the origin. Sinc $\mu$ is 1-1, part (ii) of Proposition \ref{moment} implies that we have at most $n$ of $H_k$'s.  Now the result follows from Theorems \ref{1} and \ref{3}.\hfill $\Box$

{\em Proof of Theorem \ref{6}.} Let $C^\ast(N)$ denote  the punctured metric cone over $N$. This is a hyperk\"ahler manifold which is complete except for the puncture. The moment map always exists on  a $3$-Sasakian manifold \cite{BGM} and so it exists on $C^\ast(N)$ (where it commutes with dilatations). The flats $H_k$ corresponding to $C^\ast(N)$ satisfy the condition (ii) of Proposition \ref{moment} everywhere except the origin. We can therefore, as in the proof of Theorem 1, consider moment map level sets for a generic $(n-1)$-dimensional subtorus. This shows that $C^\ast(N)$ is locally equivariantly isometric to a Taub-NUT modification of $C^\ast(N^\prime)$ where $N^\prime$ is a $3$-Sasakian quotient of a sphere determined by the $H_k$ (existence of $N^\prime$ follows from Theorem 4.1 in \cite{BD}). Such a modification, however, can be locally $3$-Sasakian if and only if it is of order zero (this follows, for example, from the fact that a higher order modification does not posess a hyperk\"ahler potential in the sense of \cite{Swann}). Therefore $C^\ast(N)$ is locally equivariantly isometric to $C^\ast(N^\prime)$.
This local isometry $f$ respects the hyperk\"ahler structures and commutes with the action of ${\Bbb R}$ by dilatations. Hence $f$ restricted to $N$ induces the $3$-Sasakian structure on its image, which must therefore be contained in $N^\prime$. Thus we have obtained an equivariant $3$-Sasakian isometry from a connected open subset of $N$ to a connected open subset of $N^\prime$. As in the proof of Theorem 1, real-analyticity, compactness and simple-connectedness of $S$ yield the result.\hfill $\Box$

{\em Proof of Corollary \ref{7}.} This follows now from Corollary 3 in \cite{BielBull}.\hfill $\Box$

\begin{ack} This work was done during my stay at the Max-Planck-Institut f\"ur Mathematik in Bonn. It is a pleasure to thank the staff and the directors of the MPI for their hospitality and support. \end{ack}

\end{document}